\documentstyle[11pt]{article}
  \newlength{\standardunitlength}
\newtheorem{cor}{Corollary} \newtheorem{lemma}{Lemma}
\newtheorem{theorem}{Theorem} 
\newenvironment{proof}{\noindent {\sc Proof:}}{$\Box$ \vspace{2 ex}}

\begin{document}

\begin{center}
A Probabilistic Approach Toward the Finite General Linear and Unitary Groups
\end{center}

\begin{center}
By Jason Fulman
\end{center}

\begin{center}
Dartmouth College
\end{center}

\begin{center}
Jason.E.Fulman@Dartmouth.Edu
\end{center}

\begin{abstract}
	Probabilistic algorithms are applied to prove theorems about the
finite general linear and unitary groups which are typically proved by
techniques such as character theory and Moebius inversion. Among the
theorems studied are Steinberg's count of unipotent elements, Rudvalis and
Shindoda's work on the fixed space of a random matrix, and Lusztig's work
on counting nilpotent matrices of a given rank.
\end{abstract}

\section{Introduction and Background} \label{Introduction}

	In thesis work done under the guidance of Persi Diaconis, Fulman
\cite{fulthesis} defined and studied measures $M_{(u,q)}$ on the set of all
partitions of all integers. The definition uses the following standard
notation from pages 3-5 of Macdonald \cite{Mac}. Recall that $\lambda$ is
said to be a partition of $n=|\lambda|$ if $\lambda_1 \geq \lambda_2 \geq
\cdots \geq 0$ and $\sum_i \lambda_i=n$. Refer to the $\lambda_i$ as the
parts of $\lambda$. Let $m_i(\lambda)$ be the number of parts of $\lambda$
equal to $i$. Then for $q>1$ and $0<u<1$ the following formula defines a
measure:

\[ M_{(u,q)}(\lambda) = [\prod_{r=1}^{\infty} (1-\frac{u}{q^r})]
\frac{u^{|\lambda|}}{q^{2 [\sum_{h<i} h m_h(\lambda) m_i(\lambda) +
\frac{1}{2} \sum_i (i-1) m_i(\lambda)^2]}\prod_i |GL(m_i(\lambda),q)|}. \]

	For $0<u<1$ and $q$ a prime power, the measures $M_{(u,q)}$ have a
group theoretic description. For this recall (for instance from Chapter 6
of Herstein \cite{Her}) that the conjugacy classes of $GL(n,q)$ are
parameterized by rational canonical form. This form corresponds to the
following combinatorial data. To each monic non-constant irreducible
polynomial $\phi$ over $F_q$, associate a partition (perhaps the trivial
partition) $\lambda_{\phi}$ of some non-negative integer
$|\lambda_{\phi}|$. Let $m_{\phi}$ denote the degree of $\phi$. The only
restrictions necessary for this data to represent a conjugacy class are:

\begin{enumerate}

\item $|\lambda_z| = 0$
\item $\sum_{\phi} |\lambda_{\phi}| m_{\phi} = n$

\end{enumerate}

	To be explicit, and for use in Section \ref{Applications}, a
representative of the conjugacy class corresponding to the data
$\lambda_{\phi}$ may be given as follows. Define the companion matrix
$C(\phi)$ of a polynomial $\phi(z)=z^{m_{\phi}} + \alpha_{m_{\phi}-1}
z^{m_{\phi}-1} + \cdots + \alpha_1 z + \alpha_0$ to be:

\[ \left( \begin{array}{c c c c c}
		0 & 1 & 0 & \cdots & 0 \\
		0 & 0 & 1 & \cdots & 0 \\
		\cdots & \cdots & \cdots & \cdots & \cdots \\
		0 & 0 & 0 & \cdots & 1 \\
		-\alpha_0 & -\alpha_1 & \cdots & \cdots & -\alpha_{m_{\phi}-1}
	  \end{array} \right) \]
	
	Let $\phi_1,\cdots,\phi_k$ be the polynomials such that
$|\lambda_{\phi_i}|>0$. Denote the parts of $\lambda_{\phi_i}$ by
$\lambda_{\phi_i,1} \geq \lambda_{\phi_i,2} \geq \cdots $. Then a matrix
corresponding to the above conjugacy class data is:

\[ \left( \begin{array}{c c c c}
		R_1 & 0 & 0 &0 \\
		0 & R_2 & 0 & 0\\
		\cdots & \cdots & \cdots & \cdots \\
		0 & 0 & 0 & R_k
	  \end{array} \right) \]

	where $R_i$ is the matrix:

\[ \left( \begin{array}{c c c}
		C(\phi_i^{\lambda_{\phi_i,1}}) & 0 & 0  \\
		0 & C(\phi_i^{\lambda_{\phi_i,2}}) & 0 \\
		0 & 0 & \cdots 
	  \end{array} \right) \]

	Now consider the following procedure for putting a measure on the
set of all partitions of all integers. Fix $u$ such that $0<u<1$. Pick a
non-negative integer such that the chance of choosing $n$ is equal to $(1-u)u^n$.
Then pick
$\alpha$ uniformly in $GL(n,q)$ and take $\lambda$ to be the paritition
corresponding to the polynomial $z-1$ in the rational canonical form of
$\alpha$ (if $n=0$ take $\lambda$ to be the trivial partition). Fulman
\cite{fulmacdonald} proves that the random partition so defined obeys
$M_{(u,q)}$ measure.

	The measures $M_{(u,q)}$ have further remarkable properties. To
state them, we use the standard notation that $\lambda_i' =
m_i(\lambda)+m_{i+1}(\lambda)+\cdots$. Then the following three equations
hold (for the third equation assume that $k \geq 2$)

\[ \sum_{\lambda} x^{|\lambda|} M_{(u,q)}(\lambda) = \prod_{r=1}^{\infty}
\frac{(1-\frac{u}{q^r})}{(1-\frac{ux}{q^r})} \]

\[ \sum_{\lambda: \lambda_1'=k} x^{|\lambda|} M_{(u,q)}(\lambda) =
\frac{(ux)^k}{|GL(k,q)|} \frac{\prod_{r=1}^{\infty}
(1-\frac{u}{q^r})}{\prod_{r=1}^k (1-\frac{ux}{q^r})} \]

\[ \sum_{\lambda: \lambda_1<k} M_{(1,q)}(\lambda) = \prod_{r=1 \atop r=0,
\pm k (mod \ 2k+1)}^{\infty} (1-\frac{1}{q^r}) \]

	These equations alone, independent of their group-theoretic
motivation, are strong evidence that the measures $M_{(u,q)}$ are worthy of
study. The first two equations will be proved probabilistically and
interpreted group theoretically in this paper. The third equation is
related to the Rogers-Ramanujan identities and has group theoretic
meaning. For more on this third equation, which seems to be the first
appearance of the Rogers-Ramanujan identities in finite group theory, see
the paper by Fulman \cite{fulRogers}.

	The measures $M_{(u,q)}$ are intimately related to symmetric
function theory, which is the topic of the companion paper by Fulman
\cite{fulmacdonald}. That paper exploits the link with symmetric functions
to develop a probabilistic algorithm, called the ``Young Tableau
Algorithm'', for growing partitions according to the measure
$M_{(u,q)}$. This and related algorithms are described in Section
\ref{Algorithms} of this paper and will be key tools of this paper.

	Section \ref{Applications} applies the tools of Section
\ref{Algorithms} to three settings. First, Steinberg's count of unipotent
elements is proved for $GL(n,q)$ and $U(n,q)$. Second, a deeper
understanding is given to work of Rudvalis and Shinoda \cite{Rud} on the
fixed space of a random element of $GL(n,q)$ or $U(n,q)$. Third, Lusztig's
\cite{Lusztig} results on nilpotent matrices of a given rank are derived
probabilistically.

	In principle, the algorithms of Section \ref{Algorithms} should be
useful for analyzing any conjugacy class function of $GL(n,q)$ or
$U(n,q)$. The results of this paper only scratch the surface.

\section{Probabilistic Algorithms} \label{Algorithms}

	Fulman \cite{fulthesis} develops probabilistic algorithms for
growing partitions according to the measures $M_{(u,q)}$. We first describe
the ``Young Tableau Algorithm''. To state it, recall that one defines the
diagram associated to $\lambda$ as the set of points $(i,j) \in Z^2$ such
that $1 \leq j \leq \lambda_i$. We follow Macdonald's convention that the
row index $i$ increases as one goes downward and the column index $j$
increases as one goes across. So the diagram of the partition $(432)$ is:

\[ \begin{array}{c c c c}
		. & . & . & .   \\
		. & . & . &     \\
		. & . &  &   
	  \end{array}  \]

	Now we state the Young Tableau Algorithm. Fulman \cite{fulthesis}
proves that for $0<u<1$ and $q>1$, the algorithm always halts and that the
resulting partition obeys $M_{(u,q)}$ measure.

\begin{center}
The Young Tableau Algorithm
\end{center}

\begin{description}

\item [Step 0] Start with $N=1$ and $\lambda$ the empty
partition. Also start with a collection of coins indexed by the
natural numbers, such that coin $i$ has probability $\frac{u}{q^i}$ of
heads and probability $1-\frac{u}{q^i}$ of tails.

\item [Step 1] Flip coin $N$.

\item [Step 2a] If coin $N$ comes up tails, leave $\lambda$ unchanged,
set $N=N+1$ and go to Step 1.

\item [Step 2b] If coin $N$ comes up heads, choose an integer $S>0$
according to the following rule. Set $S=1$ with probability $\frac
{q^{N-\lambda_1'}-1} {q^N-1}$. Set $S=s>1$ with probability
$\frac{q^{N-\lambda_s'}-q^{N-\lambda_{s-1}'}}{q^N-1}$. Then increase
the size of column $s$ of $\lambda$ by 1 and go to Step 1.

\end{description}

	The following remarks may be helpful.

\begin{enumerate}

\item Recall that a standard Young tableau $T$ of size $n$ is a partition
of $n$ with each dot replaced by one of $\{1,\cdots,n\}$ such that each of
$\{1,\cdots,n\}$ appears exactly once and the numbers increase in each row
and column of $T$. For instance,

\[ \begin{array}{c c c c c}
		1 & 3 & 5 & 6 &   \\
		2 & 4 & 7 &  &    \\
		8 & 9 &  &  &    
	  \end{array}  \]

	is a standard Young tableau. The Young Tableau Algorithm is so
named because numbering the dots in the order in which they are created
gives a standard Young tableau.

\item As an example of the Young Tableau Algorithm, suppose we are at Step
1 with $\lambda$ equal to the following partition:
	
\[ \begin{array}{c c c c}
		. & . & . & .  \\
		. & . &  &      \\
		. &  &  &    \\
		 & & &  
	  \end{array}  \]

	Suppose also that $N=4$ and that coin 4 had already come up
heads once, at which time we added to column 1, giving $\lambda$. Now
we flip coin 4 again and get heads, going to Step 2b. We add to column
$1$ with probability $\frac{q-1}{q^4-1}$, to column $2$ with
probability $\frac{q^2-q}{q^4-1}$, to column $3$ with probability
$\frac{q^3-q^2} {q^4-1}$, to column $4$ with probability $0$, and to
column $5$ with probability $\frac{q^4-q^3}{q^4-1}$. We then return to
Step 1.

\end{enumerate}

	There is a second way to decompose the measures $M_{(u,q)}$. This
uses the so-called Young lattice, which is important in combinatorics and
representation theory. The elements of this lattice are all partitions of
all numbers. An edge is drawn between partitions $\lambda$ and $\Lambda$ if
$\Lambda$ is obtained from $\lambda$ by adding one dot.

\begin{theorem} (Fulman \cite{fulthesis}) \label{Weight} Put weights
$m_{\lambda,\Lambda}$ on the Young lattice according to the rules:

\begin{enumerate}

\item $m_{\lambda,\Lambda} = \frac{u}{q^{\lambda_1'}(q^{\lambda_1'+1}-1)}$ if
$\Lambda$ is obtained from $\lambda$ by adding a box to column 1

\item $m_{\lambda,\Lambda} = \frac{u(q^{-\lambda_s'}-q^{-
\lambda_{s-1}'})}{q^{\lambda_1'}-1}$ if $\Lambda$ is obtained from
$\lambda$ by adding a box to column $s>1$

\end{enumerate}

	Then the following formula holds:

\[ M_{(u,q)}(\lambda) = [\prod_{r=1}^{\infty} (1-\frac{u}{q^r})]
\sum_{\gamma} \prod_{i=0}^{|\lambda|-1} m_{\gamma_i,\gamma_{i+1}} \]

	where the sum is over all paths $\gamma$ from the empty partition
to $\lambda$, and the $\gamma_i$ are the partitions along the path $\gamma$.
\end{theorem}

{\bf Remarks}
\begin{enumerate}

\item Observe that a Young tableau $T$ of shape $\lambda$ is equivalent to
a path in the Young lattice from the empty partition to $\lambda$. This
equivalence is given by growing the partition $\lambda$ by adding boxes in
the order $1,\cdots,n$ in the positions determined by $T$. For instance the
tableau:

\[ \begin{array}{c c c}
		1 & 3 & 4 \\
		2 &  &    \\    
	  \end{array}  \]

	corresponds to the path:

\[ \begin{array}{c c c c c c c c c c c c c c c c c c c c}
		 & & & & . & & & & . &&&& .&. & & && .&.&.                          \\
		 \emptyset & & \rightarrow &&  && \rightarrow && . && \rightarrow && . &&& \rightarrow && .& & \\    
	  \end{array}  \]

	The proof of Theorem \ref{Weight} relies on the result that the
chance that the Young Tableau Algorithm outputs the standard tableau $T$ is
equal to:

\[ [\prod_{r=1}^{\infty} (1-\frac{u}{q^r})] \prod_{i=0}^{|\lambda|-1}
m_{\gamma_i,\gamma_{i+1}} \]

	where $\gamma$ is the path in the Young lattice which corresponds
to the tableau $T$.

\item Note that in Theorem \ref{Weight} the total weight out of the empty
partition is $\frac{u}{q-1}$ and that the total weight out of any other
partition $\lambda$ is:

\begin{eqnarray*}
\frac{u}{q^{\lambda_1'}(q^{\lambda_1'+1}-1)} + \sum_{i \geq 2}
\frac{u(q^{-\lambda_s'}-q^{- \lambda_{s-1}'})}{q^{\lambda_1'}-1} & = & \frac{u}{q^{\lambda_1'}(q^{\lambda_1'+1}-1)} + \frac{u}{q^{\lambda_1'}}\\
& = & \frac{uq}{q^{\lambda_1'+1}-1}\\
& < & 1
\end{eqnarray*}

	Since the sum of the weights out of a partition $\lambda$ to a
larger partition $\Lambda$ is less than 1, the weights can also be
viewed as transition probabilities, provided that one allows for
halting.

\item Fulman \cite{fulthesis} gives a third way of building partitions
according to the measure $M_{(u,q)}$. That algorithm is similar to the
Young Tableau Algorithm, but works by adding ``horizontal strips''. It
would be interesting to find applications for this third algorithm as
well.

\end{enumerate}

\section{Applications} \label{Applications}

	This section uses the algorithms of Section \ref{Algorithms} to
obtain results about the general linear and unitary groups. The concept of
a ``cycle index'' connects the measures $M_{(u,q)}$ with these groups. The
cycle index of the general linear groups is due to Kung \cite{Kun} and
Stong \cite{St1}, though without the idea of measures on partitions. The
cycle index for the unitary groups was found by Fulman \cite{fulthesis}.

\begin{enumerate}

\item {\bf General Linear Group Cycle Index} Let $\alpha$ be an element of
$GL(n,q)$ and $\phi$ a monic, degree $m_{\phi}$ irreducible polynomial with
coefficients in $F_q$, the field of $q$ elements. Let
$\lambda_{\phi}(\alpha)$ be the partition corresponding to the $\phi$ in
the rational canonical form of $\alpha$. Then,

\[ (1-u) [1+\sum_{n=1}^{\infty} \frac{u^n}{|GL(n,q)|} \sum_{\alpha \in
GL(n,q)} \prod_{\phi \neq z} x_{\phi,\lambda_{\phi}(\alpha)}] =
\prod_{\phi \neq z} (\sum_{\lambda} x_{\phi,\lambda}
M_{(u,q^{m_{\phi}})}(\lambda)) \]

\item {\bf Unitary Group Cycle Index} The conjugacy classes of $U(n,q)
\subset GL(n,q^2)$ have a description analogous to rational canonical form
for $GL(n,q)$. Given a polynomial $\phi$ with coefficients in $F_{q^2}$ and
non-vanishing constant term, define a polynomial $\tilde{\phi}$ by:

\[ \tilde{\phi} = \frac{z^{m_{\phi}} \phi^q(\frac{1}{z})}{[\phi(0)]^q} \]

	where $\phi^q$ raises each coefficient of $\phi$ to the $q$th
power. Writing this out, a polynomial $\phi(z)=z^{m_{\phi}} +
\alpha_{m_{\phi}-1} z^{m_{\phi}-1} + \cdots + \alpha_1 z + \alpha_0$ with
$\alpha_0 \neq 0$ is sent to $\tilde{\phi}(z)= z^{m_{\phi}} +
(\frac{\alpha_1}{\alpha_0})^q z^{m_{\phi}-1}+ \cdots +
(\frac{\alpha_{m_{\phi}-1}} {\alpha_0})^qz + (\frac{1}
{\alpha_0})^q$. Fulman \cite{fulthesis} shows that all $\phi$ satisfying
$\phi = \tilde{\phi}$ have odd degree.

	Wall $\cite{Wal}$ proves that the conjugacy classes of the
unitary group correspond to the following combinatorial data. As was
the case with $GL(n,q^2)$, an element $\alpha \in U(n,q)$ associates
to each monic, non-constant, irreducible polynomial $\phi$ over
$F_{q^2}$ a partition $\lambda_{\phi}$ of some non-negative integer
$|\lambda_{\phi}|$ by means of rational canonical form. The
restrictions necessary for the data $\lambda_{\phi}$ to represent a
conjugacy class are:

\begin{enumerate}
\item $|\lambda_z|=0$
\item $\lambda_{\phi}=\lambda_{\tilde{\phi}}$
\item $\sum_{\phi} |\lambda_{\phi}|m_{\phi}=n$
\end{enumerate}

	This leads to the cycle index:

\begin{eqnarray*}
(1-u) [1+\sum_{n=1}^{\infty} \frac{u^n}{|U(n,q)|} \sum_{\alpha \in
U(n,q)} \prod_{\phi \neq z} x_{\phi,\lambda_{\phi}(\alpha)} ] & = &
\prod_{\phi \neq z, \phi=\tilde{\phi}} (\sum_{\lambda} x_{\phi,\lambda}
M_{((-u)^{m_{\phi}},(-q)^{m_{\phi}})}(\lambda))\\
& & \prod_{\phi \neq
\tilde{\phi}} (\sum_{\lambda} x_{\phi,\lambda}
M_{(u^{2m_{\phi}},q^{2m_{\phi}})}(\lambda))
\end{eqnarray*}

	It is elementary to see that $M_{(-u,-q)}$ is also a measure. The
Young Tableau Algorithm can not be applied to pick from it, however,
because some of the ``probabilities'' involved would be
negative. Nevertheless, the description in terms of weights on the Young
lattice (Theorem $\ref{Weight}$) does extend by replacing $u$ and $q$ by
their negatives.

\end{enumerate}

	The following three elementary lemmas will be of use in the
applications to follow.

\begin{lemma} \label{prodgl}

\[ 1-u = \prod_{\phi \neq z} \prod_{r=1}^{\infty}
(1-\frac{u^{m_{\phi}}}{q^{m_{\phi}r}}) \]

\end{lemma}

\begin{proof}
	For all $\phi$, perform the following substitutions in the cycle
index of the general linear groups. If $|\lambda|>0$, set
$x_{\phi,\lambda}=0$. If $|\lambda|=0$, set $x_{\phi, \lambda}=1$.
\end{proof}

\begin{lemma} \label{produn}

\[ 1-u = [\prod_{\phi \neq z, \phi=\tilde{\phi}} \prod_{r=1}^{\infty}
(1+(-1)^r \frac{u^{m_{\phi}}}{q^{m_{\phi}r}})] [\prod_{\phi \neq \tilde{\phi}}
\prod_{r=1}^{\infty} (1-\frac{u^{2m_{\phi}}}{q^{2m_{\phi}r}})] \]

\end{lemma}

\begin{proof}
	Make the same substitutions as in Lemma \ref{prodgl}, but for the
unitary groups.
\end{proof}

	Lemma \ref{HW} is elementary and is taken from page 280 of Hardy
and Wright \cite{Hardy}.

\begin{lemma} \label{HW} For real $a,y$ such that $|a| \leq 1, |y|<1$,

\[ \frac{1}{(1-ay)\cdots(1-ay^k)} = 1 + ay \frac{1-y^k}{1-y} + a^2y^2
\frac{(1-y^k)(1-y^{k+1})}{(1-y)(1-y^2)} + \cdots \]

\end{lemma}

\begin{center}
{\bf Application 1: Counting unipotent elements}
\end{center}
	
	The following theorem of Steinberg is normally proven using the
Steinberg character, as on page 156 of Humphreys $\cite{Hum}$. Recall that
$\alpha \in GL(n,q)$ is called unipotent if all of its eigenvalues are
equal to one.

\begin{theorem} \label{Steinberg} The number of unipotent elements in a
finite group of Lie type $G^F$ is the square of the order of a $p$-Sylow of
$G^F$, where $p$ is the prime used in the construction of $G^F$ (in the
case of the classical groups, $p$ is the characteristic of $F_q$).
\end{theorem}
	
	The goal of this application is to give a probabilistic proof of
Steinberg's result for $GL(n,q)$ and $U(n,q)$. To this end, we obtain a
generating function for the size of a partition $\lambda$ chosen from the
measure $M_{(u,q)}$.

\begin{theorem} \label{sizegen}

\[ \sum_{\lambda} x^{|\lambda|} M_{(u,q)}(\lambda) = \prod_{r=1}^{\infty}
\frac{(1-\frac{u}{q^r})}{(1-\frac{ux}{q^r})} \]

\end{theorem}

\begin{proof}
	Observe from the Young Tableau Algorithm that the size of the
partition is equal to the total number of coins which come up heads. The
$r=i$ term of the product on the right hand side corresponds to the tosses
of coin $i$, and these terms are multiplied because the coin tosses of
different coins are independent.
\end{proof}

\begin{cor} \label{unip} The number of unipotent elements of $GL(n,q)$ is
$q^{n(n-1)}$.
\end{cor}

\begin{proof}
	Setting $x_{z-1,\lambda}=x^{|\lambda|}$ and $x_{\phi,\lambda}=0$
for $\phi \neq z-1$ in the cycle index for $GL(n,q)$ shows that the number
of unipotent elements of $GL(n,q)$ is

\begin{eqnarray*}
& & |GL(n,q)| [u^nx^n] \frac{1}{1-u} (\sum_{\lambda} x^{|\lambda|}
M_{(u,q)}(\lambda)) (\prod_{\phi \neq z-1} \prod_{r=1}^{\infty}
(1-\frac{u^{m_{\phi}}}{q^{m_{\phi}r}}))\\
& = & |GL(n,q)| [u^nx^n] \frac{1}{1-u} (\prod_{r=1}^{\infty}
\frac{(1-\frac{u}{q^r})} {(1-\frac{ux}{q^r})}) (\prod_{\phi \neq z-1}
\prod_{r=1}^{\infty} (1-\frac{u^{m_{\phi}}}{q^{m_{\phi}r}}))\\
& = & |GL(n,q)| [u^nx^n] \prod_{r=1}^{\infty} (\frac{1}{1-\frac{ux}{q^r}})\\
& = & q^{n(n-1)}
\end{eqnarray*}

	The first equality comes from Theorem \ref{sizegen}. The second
equality is Lemma \ref{prodgl}, and the third equality is Lemma \ref{HW}
with $a=ux, y=\frac{1}{q}$.
\end{proof}

	A similar argument works for the unitary groups. It is well known
that the order of $U(n,q)$ is $q^{{n \choose 2}} \prod_{i=1}^n (q^i -
(-1)^i)$ and its $p$-Sylows have size $q^{{n \choose 2}}$.

\begin{cor} \label{unipun} The number of unipotent elements of $U(n,q)$ is
$q^{n(n-1)}$.
\end{cor}

\begin{proof}
	Setting $x_{z-1,\lambda}=x^{|\lambda|}$ and $x_{\phi,\lambda}=0$
for $\phi \neq z-1$ in the cycle index for $U(n,q)$ shows that the number
of unipotent elements of $U(n,q)$ is

\begin{eqnarray*}
& & |U(n,q)| [u^nx^n] \frac{1}{1-u} (\sum_{\lambda} x^{|\lambda|}
M_{(-u,-q)}(\lambda)) (\prod_{\phi \neq z-1, \phi= \tilde{\phi}}
\prod_{r=1}^{\infty} (1+(-1)^r \frac{u^{m_{\phi}}}{q^{m_{\phi}r}})) (\prod_{\phi
\neq
\tilde{\phi}} \prod_{r=1}^{\infty} (1- \frac{u^{2m_{\phi}}}{q^{2m_{\phi}r}}))\\
& = & |U(n,q)| [u^nx^n] \frac{1}{1-u} (\prod_{r=1}^{\infty}
\frac{(1+(-1)^r\frac{u}{q^r})} {(1+(-1)^r\frac{ux}{q^r})}) (\prod_{\phi
\neq z-1, \phi= \tilde{\phi}} \prod_{r=1}^{\infty} (1+(-1)^r
\frac{u^{m_{\phi}}}{q^{m_{\phi}r}})) (\prod_{\phi \neq \tilde{\phi}}
\prod_{r=1}^{\infty} (1- \frac{u^{2m_{\phi}}}{q^{2m_{\phi}r}}))\\
& = & |U(n,q)| [u^nx^n] \prod_{r=1}^{\infty} (\frac{1}
{(1+(-1)^r\frac{ux}{q^r})})\\
& = & q^{n(n-1)}
\end{eqnarray*}

	The first equality comes from Theorem \ref{sizegen} with $u$ and
$q$ replaced by their negatives (this substitution is valid, even though
the probabilistic proof of Theorem \ref{sizegen} breaks down for this
case). The second equality is Lemma \ref{produn}. The third equality uses
Lemma \ref{HW} with $a=-ux, y=-\frac{1}{q}$.
\end{proof}

{\bf Remarks}
\begin{enumerate}

\item The proof technique of Corollaries \ref{unip} and \ref{unipun} can be
used to give formulas for the number of elements of $GL(n,q)$ and $U(n,q)$
with a given characteristic polynomial. See Fulman \cite{fulcycle} for
details. 

\item At least for $GL(n,q)$, there should be a proof of Steinberg's count
of unipotents which is bijective (i.e. which maps a pair of elements in a
$p$-Sylow to a unipotent element).

\end{enumerate}

\begin{center}
{\bf Application 2: Work of Rudvalis and Shinoda}
\end{center}

	Rudvalis and Shinoda \cite{Rud} studied the distribution of fixed
vectors for the classical groups over finite fields. Let $G=G(n)$ be a
classical group (i.e. one of $GL$,$U$,$Sp$, or $O$) acting on an $n$
dimensional vector space $V$ over a finite field $F_q$ (in the unitary case
$F_{q^2}$) in its natural way. Let $P_{G,n}(k,q)$ be the chance that an
element of $G$ fixes a $k$ dimensional subspace and let $P_{G,\infty}(k,q)$
be the $n \rightarrow \infty$ limit of $P_{G,n}(k,q)$. Rudvalis and Shinoda
$\cite{Rud}$ obtained the following results for the general linear and
unitary cases.

\begin{enumerate}

\item $P_{GL,n}(k,q) = \frac{1}{|GL(k,q)|} \sum_{i=0}^{n-k}
\frac{(-1)^i q^{{i \choose 2}}}{q^{ki}|GL(i,q)|}$

\item $P_{GL,\infty}(k,q) =  [\prod_{r=1}^{\infty} (1-\frac{1}{q^r})]
\frac{(\frac{1}{q})^{k^2}}{(1-\frac{1}{q})^2 \cdots
(1-\frac{1}{q^k})^2} $

\item $P_{U,n}(k,q) = \frac{1}{|U(k,q)|} \sum_{i=0}^{n-k}
\frac{(-1)^i (-q)^{{i \choose 2}}}{(-q)^{ki} |U(i,q)|}$

\item $P_{U,\infty}(k,q) = [\prod_{r=0}^{\infty} (1+\frac{1}{(-q)^r})]
\frac{(\frac{1}{q})^{k^2}} {(1-\frac{1}{q^2}) \cdots (1-\frac{1}{q^{2k}})}
$

\end{enumerate}

	At first glance it is not even clear that $P_{GL,\infty}(k,q)$ and
$P_{U,\infty}(k,q)$ define probability distributions in $k$, but as
Rudvalis and Shinoda note, this follows from identities of Euler. Proofs of the
above results of Rudvalis and Shinoda used Moebius version on the lattice of
subspaces of a vector space and a detailed knowledge of geometry over finite
fields.

	Theorem \ref{Weight} of the previous section will lead to
probabilistic proofs of the above four equations. In particular, a
probabilistic interpretation will be given to the products in the formulas
for $P_{GL,\infty}(k,q)$ and $P_{U,\infty}(k,q)$. The first step is to
connect the theorems of Rudvalis and Shinoda with the partitions in the
rational canonical form of $\alpha$.

\begin{lemma} \label{RatParts} The dimension of the fixed space of an
element $\alpha$ of $GL(n,q)$ is equal to $\lambda_{z-1}(\alpha)_1'$
(i.e. the number of parts of the partition corresponding to the polynomial
$z-1$ in the rational canonical form of $\alpha$).
\end{lemma}

\begin{proof}
	It must be shown that the kernel of $\alpha-I$, where $I$ is the
identity map, has dimension $\lambda_{z-1}(\alpha)_1'$. By the explicit
description of the rational canonical form of a matrix in Section
$\ref{Introduction}$, it is enough to prove that the kernel of the linear
map with matrix $M=C((z-1)^i)-I$ is 1 dimensional for all $i$ (as in
Section $\ref{Introduction}$, $C(\phi)$ is the companion matrix of a
polynomial $\phi)$.

	Each of the first $i-1$ rows of $M$ sums to 0, and they are
linearly independent. So it needs to be shown that the last row of $M$
has sum 0. This follows from the fact that the coefficients of
$(z-1)^i$ sum to 0.
\end{proof}

	To proceed further, we need some more notation. Let $T$ be a
standard Young tableau with $k$ parts. Define numbers
$h_1(T),\cdots,h_k(T)$ associated with $T$. Let
$h_m(T)=T_{(m+1,1)}-T_{(m,1)}-1$ for $1 \leq m \leq k-1$ and let
$h_k(T)=|T|-T_{(k,1)}$. So if $k=3$ and $T$ is the tableau

\[ \begin{array}{c c c c c}
		1 & 3 & 5 & 6 &   \\
		2 & 4 & 7 &  &    \\
		8 & 9 &  &  &    
	  \end{array}  \]

	then $h_1(T)=2-1-1=0$, $h_2(T)=8-2-1=5$, and $h_3(T)=9-8=1$. View
$T$ as being created by the Young Tableau Algorithm. Then for $1 \leq m
\leq k-1$, $h_m(T)$ is the number of dots added to $T$ after it becomes a
tableau with $m$ parts and before it becomes a tableau with $m+1$
parts. $h_k(T)$ is the number of dots added to $T$ after it becomes a
tableau with $k$ parts. The proof of Theorem $\ref{Interp}$ will show that
if one conditions $T$ chosen from the measure $M_{(u,q)}$ on having $k$
parts, then the random variables $h_1(T),\cdots,h_k(T)$ are independent
geometrics with parameters $\frac{u}{q}, \cdots, \frac{u}{q^k}$. This will
explain the factorization on the right-hand side of the formula in Theorem
$\ref{Interp}$.

\begin{theorem} \label{Interp}

\[ \sum_{\lambda: \lambda_1'=k} x^{|\lambda|} M_{(u,q)}(\lambda) =
\frac{(ux)^k}{|GL(k,q)|} \frac{\prod_{r=1}^{\infty}
(1-\frac{u}{q^r})}{\prod_{r=1}^k (1-\frac{ux}{q^r})} \]

\end{theorem}

\begin{proof}
	We sum over all Young tableaux $T$ with $k$ parts $"x^{|T|}$ times
the chance that the Tableau algorithm outputs $T$". The point is that one
can easily compute the probability that the Tableau algorithm produces a
tableau $T$ with given values $h_1,\cdots,h_k$.

	Suppose that one takes a step up along the Young lattice from a
partition with $m$ parts. Theorem $\ref{Weight}$ implies that the weight
for adding to column $1$ is $\frac{u}{q^m(q^{m+1}-1)}$, and that the sum of
the weights for adding to any other column is $\frac{u}{q^m}$. Thus
$x^{|T|}$ times the chance that the Tableau algorithm yields a tableau with
given values $h_1,\cdots,h_k$ is:

\[ \prod_{r=1}^{\infty} (1-\frac{u}{q^r}) \frac{(xu)^k}{|GL(k,q)|}
\prod_{m=1}^k (\frac {ux}{q^m})^{h_m} \]

	Summing over all possible values of $h_m \geq 0$ gives:

\begin{eqnarray*}
\sum_{\lambda: \sum \lambda_1'=k} x^{|\lambda|} M_{(u,q)}(\lambda) & = &  \prod_{r=1}^{\infty} (1-\frac{u}{q^r}) \frac{(ux)^k}{|GL(k,q)|}
\prod_{m=1}^k [\sum_{h_m=0}^{\infty} (\frac{ux}{q^m})^{h_m}]\\
& = &  \prod_{r=1}^{\infty} (1-\frac{u}{q^r}) \frac{(ux)^k}{|GL(k,q)|}
\prod_{m=1}^k \frac{1}{(1-\frac {ux}{q^m})}\\
& = &  \frac{(ux)^k}{|GL(k,q)|} \frac{\prod_{r=1}^{\infty} (1-\frac{u}{q^r})}{\prod_{r=1}^k (1-\frac{ux}{q^r})}
\end{eqnarray*}

\end{proof}
 
	To deduce the The Rudvalis/Shinoda formulas for the general linear
and unitary groups, two further easy lemmas will be used.

\begin{lemma} \label{bign} If $f(1)<\infty$ and $f$ has a Taylor
series around 0, then

\[ lim_{n \rightarrow \infty} [u^n] \frac{f(u)}{1-u} = f(1) \]

\end{lemma}

\begin{proof}
	Write the Taylor expansion $f(u) = \sum_{n=0}^{\infty} a_n
u^n$. Then observe that $[u^n] \frac{f(u)}{1-u} = \sum_{i=0}^n a_i$.
\end{proof}

\begin{lemma} \label{Stong} (Goldman and Rota \cite{God})

\[ \prod_{r=1}^{\infty} (1-\frac{u}{q^r}) = \sum_{i=0}^{\infty}
\frac{(-u)^i}{(q^i-1) \cdots (q-1)} \]

\end{lemma}

	The goal of our second application of the algorithms of Section
\ref{Algorithms} can now be attained.

\begin{theorem} \label{Ru} (Rudvalis and Shinoda \cite{Rud})

\begin{enumerate}

\item $P_{GL,n}(k,q) = \frac{1}{|GL(k,q)|} \sum_{i=0}^{n-k}
\frac{(-1)^i q^{{i \choose 2}}}{q^{ki}|GL(i,q)|}$

\item $P_{GL,\infty}(k,q) =  [\prod_{r=1}^{\infty} (1-\frac{1}{q^r})]
\frac{(\frac{1}{q})^{k^2}}{(1-\frac{1}{q})^2 \cdots
(1-\frac{1}{q^k})^2} $

\item $P_{U,n}(k,q) = \frac{1}{|U(k,q)|} \sum_{i=0}^{n-k}
\frac{(-1)^i (-q)^{{i \choose 2}}}{(-q)^{ki} |U(i,q)|}$

\item $P_{U,\infty}(k,q) = [\prod_{r=0}^{\infty} (1+\frac{1}{(-q)^r})]
\frac{(\frac{1}{q})^{k^2}} {(1-\frac{1}{q^2}) \cdots (1-\frac{1}{q^{2k}})}
$

\end{enumerate}

\end{theorem}

\begin{proof}
	In the cycle index for the general linear groups, set
$x_{z-1,\lambda}=1$ if $\lambda$ has $k$ parts and $x_{\phi,\lambda}=0$
otherwise. By Lemma \ref{RatParts}, Theorem $\ref{Interp}$ with $x=1$, and
Lemma $\ref{Stong}$,

\begin{eqnarray*}
P_{GL,n}(k,q) & = & [u^n] \frac{1}{1-u} \sum_{\lambda:\lambda_1'=k}
M_{(u,q)}(\lambda)\\
& = & [u^n] \frac{u^k \prod_{r=1}^{\infty} (1-\frac{u}{q^{k+r}})}{(1-u) |GL(k,q)|}\\
& = & \frac{1}{|GL(k,q)|} [u^{n-k}] \frac{1}{1-u}
\sum_{i=0}^{\infty} \frac{(-1)^i(uq^{-k})^i}{(q^i-1) \cdots (q-1)}\\
& = & \frac{1}{|GL(k,q)|} \sum_{i=0}^{n-k} \frac{(-1)^iq^{-ki}}{(q^i-1)\cdots (q-1)}
\end{eqnarray*}

	For the second part of the theorem use Lemma \ref{bign} and Theorem
$\ref{Interp}$ with $x=1,u=1$ to conclude that:

\begin{eqnarray*}
P_{GL,\infty}(k,q) & = & lim_{n \rightarrow \infty} [u^n] \frac{1}{1-u}
\sum_{\lambda:\lambda_1'=k} M_{(u,q)}(\lambda)\\
& = & \sum_{\lambda:\lambda_1'=k} M_{(1,q)}(\lambda)\\
& = & \frac{\prod_{r=k+1}^{\infty} (1-\frac{1}{q^r})}{|GL(k,q)|}\\
& = & [\prod_{r=1}^{\infty} (1-\frac{1}{q^r})] \frac{(\frac{1}{q})^{k^2}}{(1-\frac{1}{q})^2 \cdots
(1-\frac{1}{q^k})^2}
\end{eqnarray*} 

	For the third statement, set $x_{z-1,\lambda}=1$ if $\lambda$ has
$k$ parts and $x_{\phi,\lambda}=0$ otherwise in the cycle index of the
unitary groups. As for the general linear groups,

\begin{eqnarray*}
P_{U,n}(k,q) & = & [u^n] \frac{1}{1-u} \sum_{\lambda:\lambda_1'=k} M_{(-u,-q)}(\lambda)\\
& = & [u^n] \frac{(-u)^k \prod_{r=1}^{\infty} (1-\frac{-u}{(-q)^{k+r}})}{(1-u) |GL(k,-q)|}\\
& = & \frac{1}{|U(k,q)|} [u^{n-k}] \frac{1}{1-u}
\sum_{i=0}^{\infty} \frac{(-1)^i(-u(-q)^{-k})^i}{((-q)^i-1) \cdots (-q-1)}\\
& = & \frac{1}{|U(k,q)|} \sum_{i=0}^{n-k} \frac{(-1)^i (-q)^{{i \choose 2}}}{(-q)^{ki} |U(i,q)|}
\end{eqnarray*}

	For the fourth statement argue as for the general linear groups to
conclude that:

\begin{eqnarray*}
P_{U,\infty}(k,q) & = & lim_{n \rightarrow \infty} [u^n] \frac{1}{1-u}
\sum_{\lambda:\lambda_1'=k} M_{(-u,-q)}(\lambda)\\
& = & \sum_{\lambda:\lambda_1'=k} P_{\frac{1}{(-q)^i},\frac{1}{(-q)^{i-1}},0,\frac{1}{-q}}(\lambda)\\
& = & \frac{\prod_{r=k+1}^{\infty} (1+\frac{1}{(-q)^r})}{|U(k,q)|}\\
& = & [\prod_{r=1}^{\infty} (1-\frac{1}{(-q)^r})] \frac{(\frac{1}{q})^{k^2}}{\prod_{s=1}^k (1+\frac{1}{(-q)^s}) (1-\frac{1}{(-q)^s})}\\
& = & [\prod_{r=1}^{\infty} (1+\frac{1}{(-q)^r})] \frac{(\frac{1}{q})^{k^2}}{(1-\frac{1}{q^2}) \cdots (1-\frac{1}{q^{2k}})}
\end{eqnarray*} 

\end{proof}

{\bf Remark} The reason that the formulas for $P_{GL,\infty}(k,q)$ and
$P_{U,\infty}(k,q)$ factor is because the random variables $h_i(T)$ defined
from the Young Tableau Algorithm (see the comments preceding Theorem
$\ref{Interp}$) are independent geometrics.

	Rudvalis and Shinoda obtained, after many pages of labor, the
following analogous formulas for the symplectic and orthogonal groups:

\[ P_{Sp,\infty}(k,q) = [\prod_{r=1}^{\infty} \frac{1}{1+x^r}]
\frac{x^{\frac{k^2+k}{2}}}{(1-x) \cdots (1-x^k)} \]

\[ P_{O,\infty}(k,q) = [\prod_{r=1}^{\infty} \frac{1}{1+x^r}]
\frac{x^{\frac{k^2-k}{2}}}{(1-x) \cdots (1-x^k)}. \]

	It would be marvellous if there are analogs of the Young Tableau
Algorithm for the unipotent conjugacy classes of the symplectic and
orthogonal groups. This should lead to a probabilistic interpretation of
the products in the Rudvalis/Shinoda formulas. Chapters 5 and 6 of Fulman
\cite{fulthesis} give some preliminary results in this direction.

\begin{center}
{\bf Application 3: Work of Lusztig on nilpotent matrices of
a given rank}
\end{center}

	This application uses Theorem \ref{Interp} (which was proved
probabilistically) to prove results of Lusztig \cite{Lusztig}. The two
theorems which follow were sufficiently nontrivial for Lusztig to devote a
note toward them. We found them independently.

\begin{theorem} \label{LustGL} (Lustzig \cite{Lusztig}) The number of rank
$n-k$ nilpotent $n*n$ matrices is:

\[ \frac{|GL(n,q)|}{|GL(k,q)|} \frac{(1-\frac{1}{q^k})\cdots
(1-\frac{1}{q^{n-1}})}{q^{n-k} (1-\frac{1}{q}) \cdots
(1-\frac{1}{q^{n-k}})} \]

\end{theorem}

\begin{proof}
	Adding the identity show that it suffices to count unipotent
matrices in $GL(n,q)$ with a $k$ dimensional fixed space. In the cycle
index for the general linear groups, set $x_{z-1,\lambda}=x^{|\lambda|}$ if
$\lambda$ has $k$ parts and $x_{\phi,\lambda}=0$ otherwise. By Lemma
\ref{RatParts} and Theorem $\ref{Interp}$, the sought number is:

\begin{eqnarray*}
&  & |GL(n,q)| [(ux)^n] \frac{1}{1-u} \frac{(ux)^k}{|GL(k,q)|}
\frac{\prod_{r=1}^{\infty} (1-\frac{u}{q^r})}{\prod_{r=1}^k
(1-\frac{ux}{q^r})}\\
& = & \frac{|GL(n,q)|}{|GL(k,q)|} [(ux)^{n-k}] \frac{1}{1-u}
\frac{\sum_{i=0}^{\infty} \frac{(-u)^i}{(q^i-1)\cdots (q-1)}}{\prod_{r=1}^k
(1-\frac{ux}{q^r})}\\
& = & \frac{|GL(n,q)|}{|GL(k,q)|} [(ux)^{n-k}] \frac{1}{\prod_{r=1}^k
(1-\frac{ux}{q^r})}\\
& = & \frac{|GL(n,q)|}{|GL(k,q)|} \frac{1}{q^{n-k}} \frac{(1-\frac{1}{q^k})
\cdots (1-\frac{1}{q^{n-1}})}{(1-\frac{1}{q}) \cdots (1-\frac{1}{q^{n-k}})}
\end{eqnarray*}
	
	The first equality used Lemma \ref{Stong} and the third equality
used Lemma \ref{HW}.
\end{proof}

	Theorem \ref{LustUN} is the corresponding result for the unitary
groups. Lusztig actually counted certain nilpotent matrices, but our
statement is equivalent. This can be seen using the so-called ``Cayley
Transform'' between unipotent and nilpotent matrices (page 177 of Humphreys
\cite{Hum}).

\begin{theorem} \label{LustUN} (Lustzig \cite{Lusztig}) The number of
unipotent elements of $U(n,q)$ with a $k$ dimensional fixed space is:

\[ \frac{|U(n,q)|}{|U(k,q)|} \frac{(1-\frac{1}{(-q)^k})\cdots
(1-\frac{1}{(-q)^{n-1}})}{q^{n-k} (1-\frac{1}{(-q)}) \cdots
(1-\frac{1}{(-q)^{n-k}})} \]

\end{theorem}

\begin{proof}
	Arguing as in Theorem \ref{LustGL}, the sought number is:

\begin{eqnarray*}
& & |U(n,q)| [(ux)^n] \frac{1}{1-u} \frac{(-ux)^k}{|GL(k,q)|}
\frac{\prod_{r=1}^{\infty} (1-\frac{-u}{(-q)^r})}{\prod_{r=1}^k
(1-\frac{-ux}{(-q)^r})}\\
& = & \frac{|U(n,q)|}{|U(k,q)|} [(ux)^{n-k}] \frac{1}{1-u}
\frac{\sum_{i=0}^{\infty} \frac{(u)^i}{((-q)^i-1)\cdots (-q-1)}}{\prod_{r=1}^k
(1-\frac{-ux}{(-q)^r})}\\
& = & \frac{|U(n,q)|}{|U(k,q)|} [(ux)^{n-k}] \frac{1}{\prod_{r=1}^k
(1-\frac{-ux}{(-q)^r})}\\
& = & \frac{|U(n,q)|}{|U(k,q)|} \frac{1}{q^{n-k}} \frac{(1-\frac{1}{(-q)^k})
\cdots (1-\frac{1}{(-q)^{n-1}})}{(1-\frac{1}{-q}) \cdots (1-\frac{1}{(-q)^{n-k}})}
\end{eqnarray*}

\end{proof}

\section{Acknowledgements}
	
	This work is taken from the author's Ph.D. thesis, done under the
supervision of Persi Diaconis at Harvard University. The author thanks him for
many ideas, suggestions and comments. This research was done under the generous
3-year support of the National Defense Science and Engineering Graduate Fellowship
(grant no. DAAH04-93-G-0270) and the support of the Alfred P. Sloan
Foundation Dissertation Fellowship.

\end{document}